\documentclass[12pt]{article}
\usepackage{amssymb}
\usepackage{latexsym} 
\usepackage{amsmath}
\usepackage{amscd} 
\usepackage{graphicx}
\usepackage{layout}

\newtheorem{theorem}{Theorem}[section]

\newtheorem{corollary}[theorem]{Corollary}

\newtheorem{proposition}[theorem]{Proposition}

\newtheorem{lemma}[theorem]{Lemma}

\newtheorem{rem}[theorem]{Remark}

\newcommand{\bel}{\begin{lemma}}
\newcommand{\enl}{\end{lemma}}

\newcommand{\R}{{\mathbf{R}}}

\newcommand{\N}{{\mathbf{N}}}
\newcommand{\NN}{{\mathbf{N}}}
\newcommand{\QQ}{{\mathbf{Q}}}

\newcommand{\lon}{\longrightarrow}

\newcommand{\pa}{\partial}

\newcommand{\rank}{\mathop{\rm rank}\nolimits}

\newcommand{\demo}{\par\noindent{\it Proof. \/}\ }
\newcommand{\enD}{\hfill $\Box$\vspace{3truemm} \par}

\begin{document}

\title{Generic singularities of implicit systems
of \\ 
first order differential equations on the plane}

\author{A.A. Davydov\footnote{Financial support from RFBR 00-01-00343 and
Grant-in-Aid 
for Scientific Research, No. 10304003.},
G. Ishikawa\footnote{Financial support from Grant-in-Aid 
for Scientific Research, No. 10440013.},
S. Izumiya\footnote{Financial support from Grant-in-Aid 
for Scientific Research, No. 10304003.}, 
W.-Z. Sun\footnote{Financial support from Grant-in-Aid 
for Scientific Research, No. 14604003.}
}


\renewcommand{\thefootnote}{\fnsymbol{footnote}}

\date{}

\maketitle
\begin{abstract}
For  the implicit systems of first
order ordinary differential equations on the plane there is presented
the complete local classification of generic singularities of family of its phase curves
up to
smooth orbital equivalence. Besides the well known singularities of
generic vector fields on the plane and the singularities described
by a generic first order implicit differential
equations,   
there exists only one generic singularity 
described by the implicit first order
equation supplied by Whitney umbrella surface 
generically embedded to the 
space of directions on the plane.
\end{abstract}

\section{Introduction.}

An {\it implicit system} of first order odinary differential equations
on a smooth
$n$-dimensional manifold is defined by a zero level, 
which is called {\it system surface}, 
of a smooth map $F$ from
the tangent bundle of this manifold to the $n$-dimensional Cartesian space. 
In local coordinates $x = (x_1, \dots, x_n)$ near a point of the manifold a  system
can be written in the  standard form $F(x, \dot x)=0.$

An implicit system is with
{\it locally bounded derivatives} if the restriction of the bundle
projection to the system surface  is a proper map. This restriction
is 
called 
the {\it system folding}. Only systems with locally bounded
derivatives are considered  below.
We identify the space of systems with locally bounded derivatives with
the space of respective maps $F.$ 
Then a {\it generic} system is a system
from some open everywhere dense subset in this space,  for the fine 
topology of Whitney.

A {\it solution curve} 
of an implicit system is defined as a differentiable map
$x : t \mapsto x(t)$ from an interval of the real line to the base manifold
such that the image of its natural lifting $(x(t), \dot{x}(t))$ 
to the tangent bundle belongs to the system surface. 
A {\it phase curve} is the image of such a map $x(t)$ 
and a {\it trajectory} is the image of its lifting.

In this paper we study the {\it point singularities} of the family of
phase curves provided by a germ of a system surface and present
the complete list of generic singularities on the plane 
up to smooth orbital
equivalence.

For a generic system its surface is a closed smooth $n$-dimensional submanifold
in the tangent bundle space, owing to Thom transversality theorem.
Then the system folding is a smooth map between $n$-dimensional
manifolds. Therefore the folding of a generic system can have all
generic singularities like a map between $n$-dimensional manifolds.
In fact the kernel
of the bundle projection is also $n$-dimensional 
and, due to Goryunov's theorem\cite{Go}, 
such dimension of the kernel permits all generic singularities 
for maps between  $n$-dimensional manifolds.

A generic system near a regular point of its folding can be resolved
with respect  to derivatives. In this case, near such a point, the
singularity theory 
of family of solutions of an implicit system is reduced to 
the well known theory of family of 
the phase curves 
for the generic smooth vector fields on
$n$-dimensional manifolds \cite{Ar2}.

For a generic system the velocity does not vanish at any singular point
of the system folding. Consequently, near such a point, there is well-defined
the {\it system 1-folding} which is the restriction of the projectivization
of the tangent bundle to the system surface. Again Goryunov's theorem
implies that the 1-folding of a generic system can have all
singularities like a generic map from an $n$-dimensional manifold to
an $(2n-1)$-dimensional one.

\begin{figure}[htbp]
  \begin{center}
\includegraphics[height=3truecm, width=12truecm, clip]{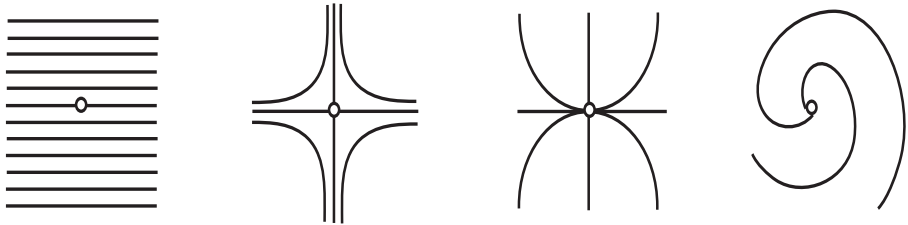} 
\caption{nonsingular point, saddle, node, and focus}
 \label{vectorfield}
\end{center}
\end{figure}

\begin{figure}[htbp]
  \begin{center}
\includegraphics[height=2truecm, width=12truecm, clip]{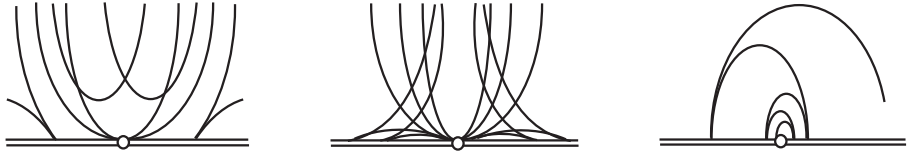} 

\ 

\caption{folded saddle, folded node, and folded focus}
 \label{foldedsing}
\end{center}
\end{figure}

\begin{figure}[htbp]
  \begin{center}
\includegraphics[height=3truecm, width=12truecm, clip]{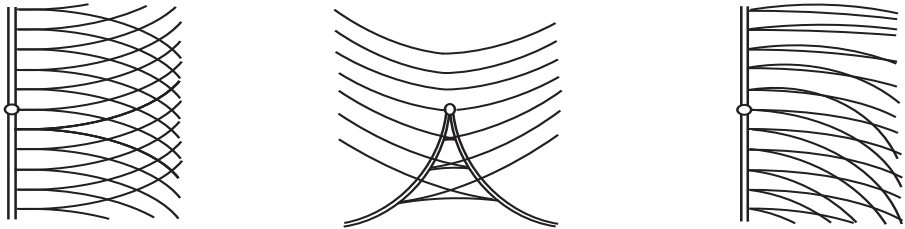} 

\ 

\caption{folded regular, pleated singular, and Whitney umbrella point}
\label{generic-fig}
\end{center}
\end{figure}

In particular, for the $2$-dimensional case the
1-folding of a generic system can have regular points and singular
points providing Whitney umbrella singularities. 
That implies, for a generic implicit system, the
classification of point singularities of families  of phase curves 
(see Theorem \ref{gen1}, \ref{gen2}). Besides well known singularities
of a generic vector fields on the plane and ones described by  generic
first order implicit differential equations there is only one
singularity provided by the implicit differential equation on the
Whitney umbrella generically embedded to the
space of directions on the plane (Figures \ref{vectorfield}, 
\ref{foldedsing}, and \ref{generic-fig}). Up to smooth orbital equivalence
the respective family of phase curve is the family of
solutions of the implicit system
$$
\dot x =\pm 1 , \qquad (\dot y)^2=x(x-y)^2
$$
near the origin. 

\ 


%

%
%
%
%
%

\ 

The family of solutions of the implicit equation
$({dy}/{dx})^2=x(x-y)^2$  has been well-known  in another situation: 
V.I.Arnol'd found it
by the investigation of slow motion of  
generic relaxation type equations
with one fast and two slow variables \cite{Ar3}\cite{Da2}. 
However the last case and the case studied in this paper are 
different. At the first place, the equation surface 
in the space of directions on the plane 
is smooth in the theory of the relaxation
type equation, while it has the Whitney umbrella singularity for the implicit
system case. 
At the second place, 
the plane distribution on the space of directions on the plane 
has singularities and is not the contact structure  in the theory of the relaxation
type equation, while it is the 
contact structure in the implicit system case. 
However if we put into correspondence to a relaxation type equation
$\dot x=\varepsilon f(x,y,z), \dot y = \varepsilon g(x,y,z), \dot z=h(x,y,z)+
\epsilon r(x,y,z)$ (where $f, g, h, r$  are smooth functions and $\varepsilon$
is a small parameter)  the surface $\dot x -f(x,y,z)=0, \dot y -g(x,y,z)=0,
h(x,y,z)=0$, then the restriction of the projection 
$(x,y,z, \dot x, \dot y) \mapsto (x,y, \dot x : \dot y)$ to this surface 
is an analogue of 1-folding. 
In a generic case this restriction is well-defined near 
any critical point of the folding being here the restriction 
of the projection $(x,y,z) \mapsto (x,y)$ to the surface $h=0$.
That reduces Arnold's case to the considered 
one. 

\ 

\begin{figure}[htbp]
  \begin{center}
\includegraphics[height=3truecm, width=12truecm, clip]{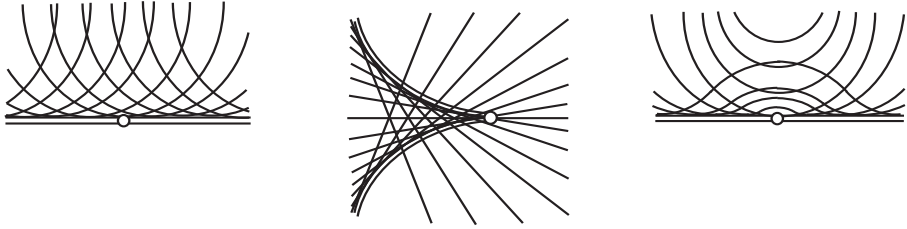} 

\ 

\caption{Clairaut fold, Clairaut cusp, and Clairaut cross cap}
\label{Clairaut-fig}
\end{center}
\end{figure}

We also present the complete classification of generic singularities
for Clairaut implicit systems on the plane (Figure \ref{Clairaut-fig}). 
An implicit system is 
called {\it of Clairaut type} if the system surface is smooth 
and,  for any critical point of its folding,  the
velocity corresponding to this point is non-zero and 
lies in the image of the tangent plane
to the system surface under the derivative of the folding, like the 
classical Clairaut equations. 
Under a mild condition, any implicit system of Clairaut type 
can be approximated by 
a system of Clairaut type 
which is foliated 
by smooth trajectories projecting to non-sigular solution 
curves. Actually we give the generic classification of such systems. 
See section \ref{Classification in the Clairaut case.} for details. 

\

%
%
%
%
%
%

\

Remark that several singularities of the systems are studied in 
the {\it parametric forms}, not in the implicit forms. 
The parametric method provides, as well as the simple process of getting 
normal forms, also the clear understanding of moduli spaces of the 
singularities and the easy  
drawing of pictures 
of the phase portaits. 

In section \ref{Classification}, we present the normal forms  both in 
generic case and generic Clairaut case. 
The proofs is given in section \ref{Generic Sing} for the generic systems. 
We give the proof in section \ref{Clairaut Sing Diag} 
for the generic Clairaut systems, 
using the theory of integral diagrams (cf. \cite{HIIY}). 

The work of this paper started during the visit of the first
author to the Department of Mathematics of Hokkaido University.
He is very thankful to the department staffs 
for a good scientific atmosphere and nice working conditions.

\section{Classification of singularities.}
\label{Classification}

Here we present the complete list of generic point singularities of first order
implicit systems on two-dimensional manifolds  
up to smooth orbital equivalence for both
general case and Clairaut one. 
Since we concern just local classification, we may treat systems on the plane 
$\R^2$.

\subsection{Classification in general case.}

\begin{theorem}
{\rm (\cite{Ar2})}
\label{gen1} 
For a generic implicit system with locally bounded derivatives 
on the plane
and for any regular point of its
folding, 
the respective point singularity takes one of the forms listed in
the second column of Table 1 
near the origin up to smooth orbital equivalence. 
Besides the parameters of the normal form from the second column  have
to satisfy the respective conditions from the third one.
\end{theorem}

\begin{rem}{\rm
In the first column of Table 1 there is pointed the standard
name for the respective singularities of  generic vector fields on
the plane.
The {\it smooth orbital equivalence} permits smooth change of coordinates and 
multiplication of vector fields by smooth positive function. 
}
\end{rem}

\begin{figure}
\begin{flushright} {\bf Table 1} \end{flushright}
$\begin{array}{|l|l|l|}
\hline
 \mbox{Type of singularities}& \mbox{Normal forms} & \mbox{Restrictions} \\
\hline
\mbox{Nonsingular point} & \dot x=1 & \\
& \dot y =0 &\\
\hline
\mbox{Nonresonance saddle} &
 \dot x=x,   & \\
\mbox{with exponent}\quad \lambda &  \dot y =\lambda y &
 \lambda \in \R_- \setminus \QQ \\
\hline
\mbox{Resonance saddle } & \dot x=x[1 \pm x^py^q +
ax^{2p}y^{2q}]& a \in \R; \, p,q \in \NN, \, 
p/q \\
\mbox{with exponent}\,-p/q & \dot y =-p y/q &
 \mbox{is an irreducible fraction. }\\
\hline
\mbox{Nonresonance node} & \dot x=\varepsilon x,
 & \\
\mbox{with exponent}\quad \lambda & \dot y =\varepsilon \lambda y&
1<\lambda \in \R_+ \setminus \N; \, \varepsilon = \pm 1 \\

\hline
\mbox{Focus with} & \dot x=\varepsilon x +\lambda y&
\lambda \in  \R_+, \varepsilon = \pm 1 \\
\mbox{exponent}\quad  \lambda & \dot y =-\lambda x +\varepsilon y & \\
\hline
\end{array}
$
\end{figure}

\begin{rem}{\rm
In \cite{Ar2} the list of the respective normal forms includes
also ones for resonance node with exponent $\lambda =n \in \NN$
( $ \dot x= x, \, \dot y =ny +\varepsilon  x^n, \,
\varepsilon \in \{-1, 0, 1\} $ ), for
resonance saddle with zero coefficients by few first resonance
monomials ($\dot x=x[1 \pm u^k + au^{2k}], \, \dot y = \lambda y, \,
1<k \in \NN $) and for degenerate  focus without formal first intergral
starting from the positive definite quadratic form ($\dot x = y \pm
x(r^{2k} +ar^{4k}), \, \dot y =  -x \pm
y(r^{2k} +ar^{4k})  $ ). All these subcases are not
generic. Ones can be removed by small perturbations of the implicit
system. }
\end{rem}

\begin{theorem}\label{gen2} For a generic implicit system on the plane
with locally bounded derivatives and any singular point of its
folding the respective point singularity takes one of the forms listed in
the second column of Table 2 
near the origin up to smooth orbital equivalence. 
Besides the parameters of the normal form from the second column  have
to satisfy the conditions from the third one.
\end{theorem}

\begin{figure}
\begin{flushright} {\bf Table 2} \end{flushright}
$\begin{array}{|l|l|l|}
\hline
 \mbox{Type of singularities}& \mbox{Normal forms} & \mbox{Restrictions} \\
\hline
\mbox{Folded regular point} & \dot x=\pm 1, \,  (\dot y)^2 =x  & \\

\hline
\mbox{Folded nonresonance } &
 \dot x=1,   &  \lambda \in \R_- \setminus \QQ;\\
\mbox{saddle with exponent}\, \lambda &  (\dot y)^2 =y-kx^2 &
 k={\lambda}/{(2\lambda+2)^2} \\
\hline
\mbox{Folded resonance saddle} & \dot x=1,\, (\dot y)^2 = y-kx^2+
& a \in \R; \, p, q \in \NN, \, 1 \neq p/q
\\
\mbox{with exponent}-p/q & +\varepsilon x
(x^{p+q}+ax^{2p+2q})&
\mbox{is an irreducible fraction};\\
&&  k=-{pq}/{(2p-2q)^2} \\
\hline
\mbox{Folded (nonresonance)} & \dot x=1,
 & \lambda \in \R_+ \setminus \NN; \\
\mbox{node with exponent}\, \lambda &  (\dot y)^2 =y-kx^2 &
 k={\lambda}/{(2\lambda+2)^2} \\

\hline
\mbox{Folded focus} & \dot x=1,&
\lambda \in  \R_+; \\
\mbox{with exponent}\, \lambda & (\dot y)^2 =y-kx^2 & k=(1+\lambda^{-2})/16 \\
\hline
\mbox{Whitney umbrella point} & \dot x=\pm 1,\,(\dot y)^2 =x(x-y)^2 &
 \\
\hline
\mbox{Pleated singular point} & \dot{x} = 1,\, x = 
\dot{y}\varphi(y, \dot{y}) 
&
\varphi \, \, \mbox{is a smooth function; } \\
&& \varphi(0,0) = \varphi_{\dot{y}}(0, 0) = 0 \\
&& \varphi_{y}(0, 0)\varphi_{\dot{y}\dot{y}}(0, 0) \not= 0 \\
\hline
\end{array}
$
\end{figure}

\begin{rem} {\rm By the topological orbital equivalence
normal forms of folded saddle, node and focus are $\dot x=1,
(\dot y)^2 = y-kx^2 $ near the origin with $k$ equals $-1, 1/20 $
and $k=1$, respectively \cite{Da1}\cite{Ku}.
}
\end{rem}

\subsection{Classification in the Clairaut case.}
\label{Classification in the Clairaut case.}

Now we denote by $T\R^2$ the tangent bundle of the plane $\R^2$, 
and by $\{ 0\}$ its zero section. We have the canonical projection 
from the tangent bundle, outside of the zero section, to 
the manifold of contact elements on the plane. 
We denote it by $\Pi : T\R^2 \setminus \{ 0\} \to PT^*\R^2$. 
Here $PT^*\R^2$ means the fiber-wise projectivization of 
the cotangent bundle $T^*\R^2$ on the plane $\R^2$. 
This projection $\Pi$ 
actually induces the $1$-folding mappings from  
surfaces in $T\R^2 \setminus \{ 0\}$. 
Remark that there does not exist a canonical isomorphism between 
the tangent bundle 
$T\R^2$ and the cotangent bundle 
$T^*\R^2$, but there exists the canonical isomorphism 
between $PT\R^2$ and $PT^*\R^2$, as the manifold of contact elements 
on the plane, defined by mapping a tangential 
direction on $\R^2$ 
to the co-direction on $\R^2$ having the direction as the kernel. 
Then the classification under the 
smooth orbital equivalence 
of surfaces in $T\R^2 \setminus \{ 0\}$ is reduced, via its $1$-folding and up to 
orientation of orbits,  
to the classification under the contact diffeomorphisms on 
$PT^*\R^2$ preserving the canonical fibration $\pi : 
PT^*\R^2 \to \R^2$. 


Let us consider a  system 
of Clairaut type in $T\R^2 \setminus \{ 0\}$ and its 1-unfolding in $PT^*\R^2$. 
Denote by $\Sigma_c$, the locus of contact singular 
points, namely, the locus on the system surface consisting 
of points where the contact form vanishes on the corresponding point in 
$PT^*\R^2$. Also denote by $\Sigma _\pi$ the locus of 
singular points on the system surface 
for the projection $\pi : T\R^2 \to \R^2$. 
Then, by definition, the system is of Clairaut type 
if and only if the folding mapping to $\R^2$ is of corank at most one 
and $\Sigma_c = \Sigma _\pi$. 

Dara \cite{Dara} gives the definition of Clairaut type equations 
for smooth surfaces in $\R^3 (\subset PT^*\R^2)$ with 
coordinates $x, y, p$, $p = \dot{y}/\dot{x}$, as follows: 
An implicit system $G(x,y,p) = 0$ is called {\it of Clairaut type in the sense 
of Dara} 
if $G_x+pG_y = AG + BG_p$ holds for some function-germs 
$A(x,y,p)$ and $B(x,y,p)$. 

By the definition, if a system is of Clairaut type 
in the sense of Dara, then we have 
$\Sigma_c = \Sigma_{\pi}$. 
In fact, in \cite{I}, it is proved that 
a non-singular system $G(x,y,p) = 0$ is of Clairaut type 
in the sense of Dara if and only if 
the system possesses the system of complete solutions 
consisting of classical (smooth) solutions. 
In particular, each trajectory projects to a non-singular 
curve via the folding. 

The converse is not true in general: 
For example, the system $G(x,y,p) = y-2p^3 = 0$ 
is not of Clairaut type 
in the sense of Dara, but it 
satisfes the condition $\Sigma_c = \Sigma _\pi$. 
Moreover it has the system of 
complete solutions $\Gamma (t,c)= (x, y, p) = (3t^2+c,2t^3,t)$ 
and each cuspidal solution curve is tangent to the discriminant $\pi(\Sigma _\pi) 
= \{ y = 0\}$ on the $(x, y)$-plane. 
Thus it is in fact a Clauraut type in the sense of this paper. 
Note that in this example the singular locus of the 
folding mapping is defined by $p^2 = 0$ on the system 
surface $\{ y = 2p^3\} = \{(x, 2p^3, p) \mid (x, p) 
\in (\R^2, 0)\}$, and the singular locus has multiple components. 

A system of Clauraut type is called {\it reduced} if the Jacobian 
of the folding mapping has no multiple components, 
or more exactly, if any differentiable function vanishing on the singular locus 
$\Sigma_{\pi}$ is devided by the Jacobian of $\pi$ restricted to the system 
surface. Then we see that any reduced system 
of Clauraut type can be approximated by a system of Clauraut type
with the property that 
each trajectory projects to a non-singular 
curve via the folding. 
Then we have 

\begin{theorem}
\label{gen3} 
A generic reduced system 
of Clairaut type on the plane
with locally bounded derivatives takes locally 
one of the forms in 
the second column of Table 3 
near the origin up to smooth orbital equivalence. 
\end{theorem}

\begin{figure}
\begin{flushright} {\bf Table 3} \end{flushright}
$\begin{array}{|l|l|l|}
\hline
 \mbox{Type of singularities} & \mbox{Normal forms} & 
\mbox{Restrictions} \\
\hline
\mbox{Nonsungular point} & \dot{x} = 1, \dot{y} = 0 &  \\
\hline
\mbox{Clairaut fold} & \dot{x} = 1, (\dot{y})^2 = y &  \\
\hline
\mbox{Clairaut cusp} & \dot{x} = 1, y = \dot{y}\varphi(x, \dot{y}) 
& \varphi(0, 0) = \varphi_{\dot{y}}(0, 0) = 0 \\
 & & 
\varphi_{\dot{y}\dot{y}}(0, 0)\varphi_{x}(0, 0) \not= 0 
\\
\hline
\mbox{Clairaut Whitney umbrella} & 
\dot{x} = 1, (\dot y)^2 = x^2y & \\
\hline
\end{array}
$
\end{figure}

The Clairaut type systems can be classified 
by considering parametric surfaces in $PT^*\R^2$. 

A {\it first order
differential equation germ (or, briefly, equation)} is defined
to be a map germ
$
f:(\R^2,0)\rightarrow J^1(\R,\R)\subset PT^*\R^2.
$
We also say that $f$ is 
{\it completely integrable} if there exists a submersion germ
$
\mu :(\R^2,0)\rightarrow \R
$
such that
$
d\mu \wedge f^*\theta=0. 
$
Here $\theta = dy - pdx$ denotes the canonical contact $1$-form on 
$J^1(\R, \R)$. 
We call
$\mu$ an {\it independent first (or, complete) integral of} $f$ and the pair
$
(\mu ,f):(\R^2,0)\rightarrow \R\times J^1(\R,\R)\subset\R\times PT^*\R^2
$
 is called {\it a holonomic system with independent first integral.}
We observe that $f|\mu ^{-1}(t)$ is 
a Legendrian immersion whose image is contained in ${\rm Image}\, f.$
If $\pi\circ f|\mu ^{-1}(t)$ are non-singular map for 
each $t\in (\R,\mu (0)),$
then $\{f|\mu ^{-1}(t)\}_{t\in\R}$ is the family of
graphs of non-singular solutions of ${\rm Image}\, f$ by the previous arguments.
We call such a system a {\it Clairaut type equation.}
These situation lead us to the following definition.  Let $(\mu ,g)$ be a pair
of a map germ $g:(\R^2,0)\rightarrow (\R^2,0)$ and a
submersion germ $\mu :(\R^2,0)\rightarrow (\R,0).$
Then the diagram
$$
(\R,0)\overset\mu\leftarrow (\R^2,0)
\overset g \rightarrow (\R^2,0)
$$
or briefly $(\mu ,g)$, is called an {\it integral diagram} if there
exists an equation $f:(\R^2,0)\rightarrow PT^*\R^2$ such that $(\mu ,f)$ is an equation germ with independent first
integral and $\pi\circ f=g,$ and we say that the integral diagram $(\mu ,g)$ is
{\it induced by} $f.$  If $f$ is a Clairaut type equation,
then $(\mu ,\pi\circ f)$ is called {\it of Clairaut type.}
Furthermore we introduce an equivalence relation among
integral diagrams.
Let $(\mu ,g)$ and $(\mu ',g')$ be integral diagrams.  Then
 $(\mu ,g)$ and $(\mu ',g')$ are {\it equivalent} ({\it respectively, strictly
equivalent}) if the diagram 
$$\CD
(\R,0) @<\mu << (\R^2,0) @>g>> (\R^2,0) \\
 @V\kappa VV         @V\psi VV        @VV\phi V \\
(\R,0) @<<\mu '< (\R^2,0)
@>>g'> (\R^2,0)
\endCD
$$
commutes for some diffeomorphism germs $\kappa,$ $\psi$
and $\phi$ (respectively, $\kappa=id_{\R}).$

\par
We give a generic classification of  Clairaut type equations 
in terms of the notion of integral diagrams which implies Theorem \ref{gen3}:

\begin{theorem} 
\label{Clairaut} 
For a generic Clairaut type equation
$$
(\mu ,f):(\R^2,0)\rightarrow \R\times J^1(\R,\R),
$$
the integral diagram $(\mu ,\pi\circ f)$ is strictly equivalent to
one of germs in the following list :
\par\noindent
{\rm (1)} $\mu = v,$ $g = (u, v)$ {\rm ; Nonsingular point.}
\par\noindent
{\rm (2)} $ \mu = v - \dfrac{1}{2}u,$ $g = (u, v^2)$ {\rm ; Clairaut fold.}
\par\noindent
{\rm (3)} $\mu _\alpha = v + \alpha\circ g\ {\rm for}\ 
\alpha \in{\mathfrak M}_{(x,y)},$
$g = (u, v^3 + uv)$ {\rm ; Clairaut cusp}.
\par\noindent
{\rm (4)} 
$\mu = v - \dfrac{1}{2}u^2,$  $g = (u, \dfrac{1}{4}v^2)$ 
{\rm ; Clairaut Whitney umbrella}.
\end{theorem}

\ 

Note that the forms (1) (2) and (4) of Theorem \ref{gen3} are obtained from those 
in Theorem \ref{Clairaut} as follows: 

(1) : The equation is given by $f = (u, v, 0)$, namely, by $p = 0$. 
Setting $\dot{x} = 1$, we have $\dot{y} = 0$. 

(2) : The equation is given by $f = (u, v^2, v)$, namely, by 
$p^2 = y$. Setting $\dot{x} = 1$, we have $(\dot{y})^2 = y$. 

(4) : The equation is given by $f = (u, \dfrac{1}{4}v^2, \dfrac{1}{2}uv)$, 
namely by $p^2 = x^2y$. Setting $\dot{x} = 1$, we have $(\dot{y})^2 = x^2y$.

The form (3) of Theorem \ref{gen3} is obtained from the parametric form 
(3) of Theorem \ref{Clairaut} as follows: 
The 
the $1$-folding map is given by 
$(u, v) \mapsto (x, y, p) = (u, v^3 + uv,  h(u, v))$, where $p = \dfrac{\dot{y}}{\dot{x}}$ 
and 
$$
h(u, v) = 
\dfrac{v - (3v^2 + u)\frac{\pa\alpha}{\pa x}(u, v^3 + uv)}{1 + (3v^2 + u)
\frac{\pa\alpha}{\pa y}(u, v^3 + uv)}. 
$$
Then $u$ and $h$ form another parametrization, and we have the implicit form 
$y = v^3 + uv = v(x, p)^3 + xv(x, p) =: \psi(x, p)$. 
By using a cusp-preserving diffeomorphism which maps 
the phase curve passing through the cusp vertex to the $x$-axis, 
we can suppose that $\psi(x, 0) = 0$. So we have 
$y = \psi(x, p) = p\varphi(x, p)$. 
By setting $\dot{x} = 1$, we have 
$y = \dot{y}\varphi(x, \dot{y})$. 

\par
Theorem \ref{Clairaut} gives a generic classification of integral diagrams of
Clairaut type under the strict equivalence.   
We remark that each germs from (1) to (4) are not equivalent
as integral diagrams. 
Thus the problem is reduced to classify germs which are contained
in the family (3) under the equivalence.
 This family is parametrized by
function germs $\alpha $ which are called {\it functional moduli.}
In \cite{IK}, 
it is given a characterization of functional moduli relative 
to the equivalence. For a map $g: x=u, y=v^3+uv$, 
define in the $x,y$-plane 
the set $\Delta$ of points where this map has three different preimages, and 
the set ${\mathcal D}$ like the boundary of the set $\Delta$. 
In \cite{Dufour} Dufour has shown the following characterization theorem.

\begin{theorem} 
Let $(\mu _\alpha ,g)$ be an integral diagram of the type 
{\rm (3)} in Theorem \ref{Clairaut}. 
Then for any $\alpha ,$ there exists a function germ $\alpha ':(\R^2,0)
\rightarrow (\R,0)$ such that
\par\noindent
{\rm (1)} $(\mu _\alpha ,g)$ is equivalent to $(\mu _{\alpha '},g).$
\par\noindent
{\rm (2)} $\alpha '|{\mathcal D}=0.$ 
\end{theorem}

\par

Dufour has also shown that the uniqueness of functional 
moduli relative to the equivalence. 
 We say that $\alpha$ and $\alpha '$ are {\it equivalent as moduli} if there exists
$a\in\R\setminus\{0\}$ such that 
$
a\alpha (x,y)=\alpha '(a^2x,a^3y)
$
for any $(x,y)\in \Delta.$
We remark that the above definition of the equivalence among functional moduli
is slightly different from Dufour's original definition in \cite{Dufour}.
If we adopt his definition, we cannot assert the necessity of the condition that
functional moduli are equivalent.
Actually, in \cite{IK}, we have introduced the above definition 
and shown the following theorem:

\begin{theorem} 
Let $(\mu _\alpha ,g)$ and $(\mu _{\alpha '},g)$ be integral diagrams
of {\rm (3)} such that $\alpha |{\mathcal D}=\alpha '|{\mathcal D}=0.$
Then $(\mu _\alpha ,g)$ and $(\mu _{\alpha '},g)$ are equivalent as integral diagram if and only if $\alpha $
and $\alpha '$ are equivalent as moduli.
\end{theorem}

\par
This theorem asserts that the equivalence classes of functional moduli $\alpha $ with
$\alpha |{\mathcal D}=0$ are the complete invariant
for generic classifications of Clairaut type equations under the equivalence relation
given by the group of point transformations.
\par
We define
$
{\mathfrak M} ({\mathcal D})=\{\alpha\in {\mathfrak M}_{(x,y)} 
\mid \alpha |{\mathcal D}=0\}
$
and ${\mathcal M}_{\rm cusp}={\mathfrak M}({\mathcal D})/\!\sim,$
where $\sim$ denotes the equivalence relation as moduli.
The above theorem asserts that the moduli space for
generic Clairaut type equation is ${\mathcal M}_{\rm cusp}.$

\begin{rem}{\rm 
In the paper \cite{HIIY}, it is given the classification 
of first order implicit differential equations in $PT^*\R^2$ 
endowed with independent first integrals. 
In \cite{HIIY}, the pleated singular points are called the {\it regular cusps} 
and their normal forms are give in the parametric forms: 
$(u, v) \mapsto (u^3 + uv, v)$, the parametrization of the folding map, 
with the first integral 
$\mu = \dfrac{3}{4}u^4 + \dfrac{1}{2}u^2v + 
\alpha(u^3 + uv, v)$, where $\alpha$ is any function with $\alpha(0, 0) = 0, 
\dfrac{\pa\alpha}{\pa y}(0, 0) = \pm 1$ (Theorem B (5) of \cite{HIIY}). 
Moreover, by a theorem of Kurokawa (Theorem A of \cite{K}), 
we can take, up to the equivalence, the functional moduli $\alpha$ 
satisfying 
$$
\alpha(0, y) = \pm y + \dfrac{1}{2}\dfrac{\pa^2 \alpha}{\pa y^2}(0, 0)y^2, (y \leq 0). 
$$
Furthermore the sign $\pm 1$ and $\chi_\alpha 
= \dfrac{\pa^2 \alpha}{\pa y^2}(0, 0)$ are invariants of the equation (cf. \cite{K}). 

Then the $1$-folding map is given by 
$(u, v) \mapsto (x, y, p) = (u^3 + uv, v, k(u, v))$, where $p = \dfrac{\dot{y}}{\dot{x}}$ 
and 
$$
k(u, v) = 
\dfrac{u + \frac{\pa\alpha}{\pa x}(u^3 + uv, v)}{\frac{1}{2}u^2 - 
\frac{\pa\alpha}{\pa y}(u^3 + uv, v)}. 
$$
Remark that functions 
$v$ and $k$ provide another parametrization of 
the system surface.  
Then we have 
$x = u^3 + uv = u(v, k)^3 + u(v, k)v = u(y, p)^3 + 
u(y, p)y =: \psi(y, p)$. 
Note that the locus $p = 0$ defines a smooth curve that is tangent to the $y$-axis 
at the cusp vertex on $(x, y)$-plane. By using a cusp-preserving diffeomorphism 
which maps this curve to the $y$-axis, we may suppose 
$\psi(y, 0) = 0$. Thus we have the form 
$x =  p\varphi(y, p)$. Then we have the normal form in table 2 by 
setting 
$\dot{x} = 1$.

The first integral  
$\mu_{\alpha} = \dfrac{3}{4}u^4 + \dfrac{1}{2}u^2v + 
\alpha(u^3 + uv, v)$ for the system with the folding map 
$(u, v) \mapsto (u^3 + uv, v)$ is transformed 
to the first integral  
$\mu_{\alpha'} = \dfrac{3}{4}u^4 + \dfrac{1}{2}u^2v + \alpha'(u^3 + uv, v)$ 
for the system with the same folding map 
$(u, v) \mapsto (u^3 + uv, v)$ if and only if 
$\alpha'(x, y) = \alpha(x, y)$ or $\alpha'(x, y) = \alpha(-x, y)$ 
on the cuspidal open domain $\Delta := \{ y^3 + \dfrac{27}{4}x^2 < 0\}$ 
near the origin (Proposition 6.1 of \cite{HIIY}). 

The exact description of the moduli space of pleated singular points 
remains open. However 
the reason of the existence of the functional moduli clearly 
comes from that 
for the $3$-webs on the plane by solution curves, and we see, as mentioned 
above,  
the moduli space is dominated by the function space 
$$
\left\{ 
\alpha
\vert
\overline{\Delta}
\ \left\vert \  \alpha : (\R^2, 0) \to (\R, 0), 
\alpha(0, y) = \pm y + \dfrac{1}{2}\dfrac{\pa^2 \alpha}{\pa y^2}(0, 0)y^2, (y \leq 0) 
\right\}
\right.
. 
$$
}
\end{rem}

\section{Singularities in a general case.}
\label{Generic Sing}

Here we prove Theorem \ref{gen1}, \ref{gen2}. At the first place we
study typical singularities
of system (1-)folding, and then on that base  the theorems are proved.

\subsection{Generic singularities of 1-folding.}

\begin{proposition}
\label{prop1}
For a generic implicit system its surface is either empty or a smooth
two dimensional submanifold of the tangent bundle space.
\end{proposition}

Proposition \ref{prop1} follows immediately from Thom transversality theorem.

\begin{proposition}\label{prop2}
The folding of a generic implicit system with locally bounded derivatives
is $LR$-stable map and it can have singular points of type
Whitney either fold or pleat only. In other words, 
near any its critical point,  the
folding takes respectively the form
$$
\mbox{either}\qquad  \left\{\begin{array}{l} x=u\\ y=v^2 \end{array} \right.
\qquad \mbox{or}\qquad  \left\{\begin{array}{l} x=u\\ y=v^3+uv
 \end{array}\right.
 $$
in appropriate smooth coordinates near this point and its image under
the folding with the origin at them.
\end{proposition}

Recall that  the system folding is the restriction of the tangent bundle
projection to the system surface. Thus due to Goryunov theorem such
restriction  in a generic case can have all generic singularities as a
generic map between $n$-dimensional manifolds
with the dimension of the kernel being
no greater then the dimension of the kernel of this projection
\cite{Go}. But the last dimension is also equal to $n$, and that
 permits all generic
singularities between  the system surface and the phase
space.

In the two dimensional case these singularities are Whitney fold and
pleat  \cite{Wh}, \cite{Ar1}. Besides our systems are with
locally bounded derivatives. Therefore
system foldings are  proper maps. Consequently, for a generic system its
folding is $LR$-stable map \cite{GG}, \cite{Ma}.

Thus Proposition \ref{prop2} is true.

\begin{rem} {\rm A map is called $LR$-{\it stable} (= {\it left right
stable}) if for any map being sufficiently close to it these two maps
can be carried one to another by diffeomorphisms of the image space
and the preimage space which are close to the identities. For example
the map $(x,y) \mapsto x$ from the circle $x^2+y^2=1$ to the $x$-axis
is $LR$-stable (under small perturbations in $C^k$-topology with $k\ge 2$.}
\end{rem}

\begin{proposition}\label{prop3}

For a generic implicit system with locally bounded derivatives
any critical point of its folding does not belong to the zero
section of the tangent bundle to the phase space.
In particular, near such a point the system 1-folding is well
defined.
\end{proposition}

This proposition follows immediately from Thom transversality theorem
because the conditions

$$F=0, \quad \rank F_{\dot x} <n, \quad \dot x =0$$
defines in the jet space of systems the closed Whitney stratified
manifold of the codimension $2n+1$ which is greater then the dimension
$2n$ of the tangent bundle. Therefore for a generic system these
conditions can not be satisfied semiltaniously.

\begin{theorem}\label{th3}
For a generic implicit system with locally bounded derivatives
any critical point of its folding is either regular point or
critical point of type Whitney umbrella
for the 1-folding of this system. Besides
the image of the set of critical points of the system folding
under the 1-folding is generically replaced with respect to the
standard contact structure in the  space of directions on the phase space.
In particular,
it can have tangency with this structure only of the first order and only
at  the points being regular for the system 1-folding and critical one
of type Whitney fold for the system folding.
\end{theorem}

\begin{rem}{\rm
A generic replacement of the image with respect to
the direction axis in the projectivization of the tungent bundle
follows  immediately from Proposition \ref{prop2} and the following
corollary of the first statement of Theorem \ref{th3}.
}
\end{rem}

\begin{corollary}\label{cor1}
For a generic implicit system with locally bounded derivatives
any critical point of type Whitney umbrella this system 1-folding
is the critical point of type Whitney fold for this system
folding.
\end{corollary}

Theorem \ref{th3} also implies immediately
\begin{corollary}\label{cor2}
For a generic implicit system with locally bounded derivatives
its point singularities provided by the critical point of this system
folding are described either by generic singularities of first order
implicit differential equation or by such an equation provided by
a generic replacement in the space of directions on the plane of germ
of the Whitney umbrella at it vertex.
\end{corollary}

Let us prove Theorem  \ref{th3}. Let the surface of a generic
implicit system with locally bounded derivatives be not empty.

For such a system its folding is $LR$-stable due to Proposition
\ref{prop2} and in the strength of Proposition \ref{prop3} the set of critical
points of this folding does not intersect the zero section of the tangent
bundle to the phase space.  Therefore for any system being
sufficiently close to the given one the 1-folding is well defined
near this set. 

Note that any regular point of this system folding is also
regular point of its 1-folding because the folding provides two
components of the 1-folding.

Now we again can apply Goryunov theorem \cite{Go}. Due to this theorem the
1-folding of a generic system can have all generic singularities
like a generic map from two-dimensional manifold to 3-dimensional one.
But any critical point of such a map is of Whitney umbrella type. Taking
into account Corollary \ref{cor1} we find that near such a point of a
generic system this system 1-folding takes form

$$
   x=u, \qquad y=v^2, \qquad z=uv
$$
in appropriate smooth coordinates near this critical point
and local smooth coordinates in the image space fibered over the phase
space ($x,y$-space) with the origins at this point and its image,
respectively.

A typical replacement (of
the image of the set of critical points of the folding
under the 1-folding with respect to the standard
contact structure in the space of directions on the phase space) can be
obtained by  small rotations of the tangent planes. Really such a
rotations provides small perturbation of the system but they do not
change the set of critical value of the system folding and can supply
any small rotation of the field of direction defined by our system on
this set. Finally if for a generic system the "replacement" is typical
then it is also typical for any system being sufficiently close to the
chosen one due to $LR$-stability of the folding of a generic system
in the strength of Proposition \ref{prop2}.

Theorem \ref{th3} is proved.

\subsection{Proofs of Theorems \ref{gen1}, \ref{gen2}.}

At the first place we prove Theorem \ref{gen1} and then Theorem \ref{gen2}.

Let $P$ be a regular point of the folding of a generic system. Near such a
point this system surface is smooth section of the tangent bundle.
This section provides smooth vector field $v$ near the image  $\tilde P$
of this point under the system folding.

If the point $P$ does not belong to the zero section of the tangent
bundle then this field does not vanish at the point $\tilde P$.
In that case the germ of the field $v$ at this point is
$C^{\infty}$-diffeomorphic to the germ of the constant vector field
$(1, 0)$ at the origin \cite{Ar2}. That gives the first singularity
from Table 1.

If the point $P$ belongs to zero section of the tangent bundle then field
$v$ vanishes at the point $\tilde P$. Small perturbations of the studied
system implies smooth small changing  of the field near the point $P$ due
to $LR$-stability of this system folding in accordance with
Proposition \ref{prop3}. But for a generic system this map
("small perturbation" $\mapsto$
"small perturbations $v$ near $\tilde P$") is continuous and small
perturbations of the system provides all small perturbations of the system
surface near the point $P$ and, hence, all small perturbations of the field
$v$ near the point $\tilde P$. Therefore for a generic system vector field
$v$ has at a point $\tilde P$ a generic singular point. Now  the rest
part of Theorem \ref{gen1} follows from the classical results about
normal forms of generic vector fields near singular points
up to smooth orbital equivalence \cite{Ar2}.
Theorem \ref{gen1} is proved.

Let us prove Theorem \ref{gen2}.
Due to Corollary \ref{cor2} local singularities of a generic implicit system
at singular points of its folding are described by generic singularities of
first order implicit ODE except the singular points of this system 1-folding
of type Whitney umbrella. But here we need to take into account that
a critical point of a generic implicit system never belongs the zero section
of the tangent bundle due to Proposition \ref{prop3}. Consequently, generic
singularities of implicit equations of type folded regular point,
folded singular point and pleated point gives here first five singularities
and the last one from Table 2.

Now due to Corollary \ref{cor2} to finish the proof one need to get the normal
form of an implicit
first order  ODE provided by a generic replacement of Whitney umbrella to
the space of directions on the plane. For a generic system such a replacement
has to have the following properties.  At the Whitney umbrella vertex both the
the contact plane and the vertical direction do not tangent
to the image under a generic system 1-folding of set of critical points of
this
system folding and the "handle" of this umbrella.

Consequently, smooth local coordinate systems $u,v$ on the system surface
and $x,y$ on the phase space with the origins at the studied point and its
image
under the system folding can be choosen such that this system 1-folding takes
the form
\begin{equation}\label{eq1}
 x=v^2, \qquad y=u, \qquad \frac{dy}{dx}=h(u,v)
\end{equation}
where $h$ is a smooth function, $h(0,0)=0$,
and ${dy}/{dx}$ is local coordinate along the direction axis and also
in these coordinates the "handle" of the Whitney umbrella is over the line
$x-y=0$.

Near the origin on the $y$-axis and the line $\{x-y=0\}
\cup \{x>0\},$ the studied implicit first order equation gives two
smooth direction fields
$dy/dx = f_1(y)$ and $dy/dx=f_2(y)$  respectively, where $f_1, f_2$
are smooth functions; $f_1(0)=f_2(0)=0$ because $h(0,0)=0$.
Due to Hadamard lemma these functions can be presented in the form
$f_1(y)=y\tilde f_1(y), f_2(y)=y\tilde f_2(y)$ where $\tilde f_1, \tilde f_2$
are smooth functions.
Near the origin the direction field
$$ dy/dx = (y-x)\tilde f_1(y) + x\tilde f_2(y)$$
provides the semiltanious smooth extension of these two fields.

Near the origin the extended field has a first integral of the form
$y +xI_1(x,y)$ where $I_1$ is smooth function. Taking this integral and
the function  $u+v^2I_1(v^2,u)$ like new coordinates $y$ and $u$,
respectively, one preserves the first two forms of the equation 
(\ref{eq1}) but in new coordinates the function $h$ is identically
zero on the $u$-axis and on the set corresponding to the "handle".
In new coordinates near the origin  the last set
can be defined by some equation $u-v^2X(v^2)=0$ where $X$ is a smooth
function, $X(0)>0$. Due to Hadamard lemma the function $h$ in
the  last equation from
(\ref{eq1}) can be written in the form $h(u,v)=v(u-v^2X(v^2))H(u,v)$
where $H$ is smooth function; $H(0,0) \neq 0$ because at the studied point
the system 1-folding has singularity of type Whitney umbrella.
Consequently near the origin the rescalings $\tilde v=v\sqrt{X(v^2)}$ and
$\tilde x = xX(x)$ reduce the system (\ref{eq1}) of equations to the form
("tilde" in notations of new coordinates is omitted)

\begin{equation}\label{eq2}
 x=v^2, \qquad y=u, \qquad \frac{dy}{dx}=v(u-v^2)H(u,v)
\end{equation}
with some new smooth function $H$ no vanishing at the origin.
It is easy to see that the direction field  provided by
the last equation can be lifted to the smooth direction field
\begin{equation}\label{eq3}
du/dv =2v^2(u-v^2)H(u,v).
\end{equation}
The last direction field has first integral of the form
$$I(u,v)=u+v^3J_1(u,v)$$
as it is easy to see.

Now it is sufficient to get normal form of this integral by the
changing of coordinates commuting with the involution $(u,v) \mapsto
(u,-v)$ defined by the system folding.
Taking new coordinate $u$ in  the form  $(I(u,v)+I(u,-v))/2$
(being even with respect to $v$ and so it is permitted by our involution)
we reduce the integral to the form
$I(u,v)=u +v^3J_1(u,v^2)$ or
$$I(u,v)=u+ v^3a(u)+v^5b(u) +v^7c(u,v^2),$$
where $a, b, c$ are some smooth functions; $a(0)=0\neq a'(0)b(0)$ because
$v^2u$
is the term of the lower degree in the right hand side of the equation
(\ref{eq3}), $b(0) \neq 0$ due to $H(0,0)\neq 0.$

The following lemma completes the proof

\begin{lemma} {\rm {(\cite{Ar3}, \cite{Da2})}} Near the origin the
$(u,v)$-plane
a function $u+ v^3a(u)+v^5b(u) +v^7c(u,v^2)$ with smooth
functions $a,\, b$ and $c$, $a(0)=0\neq a'(0)b(0)$, is reduced to the
form $u+ v^3u+v^5$ by a smooth diffeomorphism preserving the origin and
commuting with the involution $(u,v) \mapsto (u,-v)$.
\end{lemma}

Hence the function $H$ in the equation (\ref{eq3}) can be reduced to $1$. After
that this  equation takes the form $du/dv =2v^2(u-v^2)$. Near the points under
consideration that implies the equation $dy/dx=v(u-v^2)$ on the system
surface or the equation
$(dy/dx)^2=x(y-x)^2$ on the $(x,y)$-plane.
That gives the rest (fifth) singularity from the list of Table 2.
Theorem \ref{gen2} is proved.

\section{Clairaut type equations and Legendre singularity theory.} 
\label{Clairaut Sing Diag}

\subsection{Legendrian unfoldings.}

We briefly review the theory of one-parameter Legendrian unfoldings.
We now consider the $1$-jet bundle $J^1(\R\times\R,\R)$
and the canonical $1$-form $\Theta$ on the space.  Let
$(t,x)$ be the canonical coordinate on $\R\times\R$
and
$
(t,x,y,q,p)
$
be the corresponding coordinate on $J^1(\R\times\R,\R).$
Then the canonical $1$-form is given by
$
\Theta = dy-pdx-qdt=\theta -qdt.
$
We also have the natural projection
$$
\Pi :J^1(\R\times\R,\R)\rightarrow \R\times\R\times \R
$$
defined by 
$
\Pi (t,x,y,q,p)=(t,x,y).
$
We call the above $1$-jet bundle {\it an unfolded $1$-jet bundle.}
Let $(\mu ,f)$ be an equation with complete integral. Then
there exists a unique function germ $h :(\R^2,0)\lon \R$ such that
$
f^*\theta =h\cdot d\mu.
$ Define a
map germ
$$
\ell _{(\mu ,f)}:(\R^2,0)\rightarrow J^1(\R\times\R,\R)
$$
by 
$$
\ell _{(\mu ,f)}(u)=(\mu (u),x\circ f(u),y\circ f(u),h(u),p\circ f(u)).
$$
Then we can easily show that if $(\mu ,f)$ is a Clairaut type equation, then 
$\ell _{(\mu ,f)}$ is a Legendrian immersion germ.  
We call $\ell _{(\mu ,f)}$ {\it a complete 
Legendrian unfolding associated with} $(\mu ,f).$ 
By the aid of the notion of Legendrian unfoldings, 
Clairaut type equations are characterized as follows :

\begin{proposition} 
\label{immersion}
Let
$(\mu ,f):(\R^2,0)\rightarrow \R\times J^1(\R,\R)\subset PT^*\R^2$
be an equation with complete integral.  Then $(\mu ,f)$ is
a Clairaut type equation if and only
if $\ell _{(\mu ,f)}$ is a Legendrian non-singular Legendrian immersion germ.
\end{proposition}

A complete Legendrian unfolding $\ell _{(\mu ,f)}$ associated to 
$(\mu ,f)$ is called {\it a Legendrian unfolding of Clairaut type} if
$\ell _{(\mu ,f)}$ is a Clairaut type equation.

\subsection{Genericity.}

\par
Returning to the study of equations with complete integral, we now establish
the notion of the genericity. 
\par
Let $U\subset \R^2$ be an open set. We denote by 
$
{\rm Int}(U,\R\times J^1(\R,\R))
$
the set of systems with complete integral
$
(\mu ,f):U\rightarrow \R\times J^1(\R,\R).
$
We also define
$
L(U,J^1(\R\times\R,\R))
$
to be the set of complete Legendrian unfoldings
$
\ell _{(\mu ,f)}:U\rightarrow J^1(\R\times\R,\R).
$
\par
These sets are topological spaces equipped with the Whitney
$C^\infty$-topology.  A subset of 
either spaces is said to be {\it generic} if it is an open dense subset in the
space.
 \par
 The genericity of a property of germs are defined as follows.
Let $P$ be a property of equation germs with complete integral
$
(\mu ,f):(\R^2,0)\rightarrow \R\times J^1(\R,\R)
$
(respectively, Legendrian unfoldings
$
\ell _{(\mu ,f)}:(\R^2,0)\rightarrow
J^1(\R\times\R,\R)).
$
For an open set $U\subset \R^2,$ we define ${\mathcal P}(U)$ to be the set
of
$
(\mu ,f)\in \text{Int}(U,\R\times J^1(\R,\R))
$
(respectively,
$
\ell _{(\mu ,f)}\in L(U,J^1(\R\times\R,\R)))
$
such that the germ at $x$ whose representative is given by
$(\mu ,f)$ (respectively, $\ell _{(\mu ,f)}$) has property $P$ for any
$x\in U.$
\par
The property $P$ is said to be {\it generic} if for some neighbourhood
$U$ of $0$ in $\R^2,$ the set ${\mathcal P}(U)$ is a generic subset in
$\text{Int}(U,\R\times J^1(\R,\R))$ (respectively,
$L(U,J^1(\R\times\R,\R)).$
\par
By the construction, we have a well-defined continuous mapping
$$
(\Pi _1)_*:L(U,J^1(\R\times\R,\R))\rightarrow
\text{Int}(U,\R\times J^1(\R,\R))
$$
defined by
$
(\Pi _1)_*(\ell _{(\mu ,f)})=\Pi _1\circ\ell _{(\mu ,f)}=(\mu ,f),
$
where $\Pi _1:J^1(\R\times\R,\R)\rightarrow 
J^1(\R,\R)$ is the canonical projection.  Then it has been shown 
the following fundamental theorem: 

\begin{theorem}
The continuous map 
$$
 (\Pi _1)_*:L(U, J^1(\R\times\R, \R))\rightarrow
 \text{\rm Int}(U, \R\times J^1(\R,\R)) 
$$
is a homeomorphism.
\end{theorem}

As in \cite{IK}, we define 
the equivalence relation among parametric systems under
the group of point transformations: Two equations 
$f, f' : (\R^2, 0) \to PT^*\R^2$ are {\it equivalent under
the group of point transformations} 
if there exists a diffeomorphism $\phi : (\R^2, \pi(f(0))) \to 
(\R^2, \pi(f'(0)))$ such that the canonical lifting 
$\hat{\phi} : (PT^*\R^2, f(0)) \to (PT^*\R^2, f'(0))$ transforms 
$f$ to $f'$, namely, $\hat{\phi}\circ f = f'\circ\psi$ 
for some diffeomorphism $\psi : (\R^2, 0) 
\to (\R^2, 0)$. 

We need the following basic result: 

\begin{proposition}
Let $f, f' : (\R^2, 0) \to PT^*\R^2$ be completely integrable 
equations with independent first integrals $\mu, \mu' : (\R^2, 0) 
\to (\R, 0)$ respectively. Assume 
$\pi\circ f$ and $\pi\circ f' : (\R^2, 0) \to (\R^2, 0)$ 
have nowhere dense singular sets. 
Then $f$ and $f'$ are equivalent under
the group of point transformations if and only if the induced
integral diagrams $(\mu ,\pi\circ f)$ and $(\mu ',\pi\circ f')$ are
equivalent. 
\end{proposition}

\demo 
Assume $(\mu, \pi\circ f)$ and $(\mu ', \pi\circ f')$ are
equivalent by diffeomorphisms $(\kappa, \psi, \phi)$. 
Then 
$\phi$ maps integral curves $\pi(f(\mu^{-1}(\mu(u, v))))$ through  
$\pi(f(u, v))$ 
to $\pi(f'({\mu'}^{-1}(\kappa(\mu(u, v))))) = \pi(f'({\mu'}^{-1}(\mu'(\psi(u, v)))))$ 
through $\pi(f'(\psi(u, v)))$, so the tangent lines to them. 
Since the set of contact singular points is contained in the set of critical points 
of the projection $\pi$, we see $\hat{\phi}\circ f = f'\circ\psi$. 
This implies $f$ and $f'$ are equivalent under
the group of point transformations. 
The converse implication is clear. 
\enD

\section{Proofs for Theorem \ref{Clairaut} and Theorem \ref{gen3}.} 

In the case when $\ell _{(\mu ,f)}$ is a Legendrian immersion germ,
there exists a
generating family of $\ell _{(\mu ,f)}$ by the Arnol'd-Zakalyukin's theory
(\cite{Ar1}).  In this case the generating family is naturally
constructed by a one-parameter family of generating families associated with
$(\mu ,\ell ).$  Let
$
F:((\R\times\R)\times\R^k,0)\rightarrow (\R,0)
$
be a function germ such that
$d_2F|0\times\R\times\R^k$
is non-singular, where
$
d_2F(t,x,q)=(\frac{\partial F}{\partial q_1}(t,x,q), \dots ,
\frac{\partial F}{\partial q_k}(t,x,q)).
$
We call $F$ {\it a Morse family.}  Then $C(F)=d_2F^{-1}(0)$ is
a smooth surface germ and $\pi _F:(C(F),0)\rightarrow \R$ is a
submersion germ, where $\pi _F(t,x,q)=t.$  We call the submanifold
$C(F)$ {\it a catastrophe set of } $F.$  Define
$$
\tilde\Phi _F:(C(F),0)\rightarrow J^1(\R,\R)
$$
by
$$
\tilde\Phi _F(t,x,q)=(x,F(t,x,q),\frac{\partial F}{\partial x}(t,x,q))
$$
and
$$
\Phi _F:(C(F),0)\rightarrow J^1(\R\times\R,\R)
$$
by
$$
\Phi _F(t,x,q)=(t,x,F(t,x,q),\frac{\partial F}{\partial t}(t,x,q),
\frac{\partial F}{\partial x}(t,x,q)).
$$
Since $\frac{\partial F}{\partial q_i}=0$ on $C(F),$ we can easily show that
$
(\tilde\Phi _F)^*\theta = \frac{\partial F}{\partial t}|C(F)\cdot
dt|C(F)=0.
$
By definition, $\Phi _F$ is a Legendrian unfolding associated with the
Legendrian family $(\pi _F,\tilde\Phi _F).$  By the same method of the
theory of Arnol'd-Zakalyukin (\cite{Ar1}), we can show the following
proposition.

\begin{proposition} All Legendrian unfolding germs are constructed
by the above method.
\end{proposition}

 Let $(\mu ,f)$ be a Clairaut type equation. 
By Proposition \ref{immersion}, $\ell _{(\mu ,f)}$ is 
a Legendrian immersion.
Then we can choose a family of function germs 
$$
F:(\R\times\R,0)\rightarrow (\R,0)
$$
such that ${\rm Image}\, j^1F_t=f(\mu ^{-1}(t))$ for any 
$t\in \R,$
where $F_t(x)=F(t,x).$
We remark that
the map germ
$$
j^1_1F:(\R\times\R,0)\rightarrow J^1(\R,\R)
$$
defined by $j^1_1F(t,x)=j^1F_t(x)$
is not necessary an immersion germ.
In this case we have $(C(F),0)=(\R\times\R,0)$ and
$$
\Phi _F=j^1F:(\R\times\R, 0)\rightarrow
J^1(\R\times\R, \R),
$$
so that it is a complete Legendrian unfolding associated with
$(\pi _1,j^1_1F).$
Thus the generating family of a Legendrian unfolding of Clairaut type is given by the
above germ.
\par
In order to prove Theorem \ref{Clairaut}, 
we now introduce equivalence relations among Legendrian unfoldings.
 Let $(\mu ,g)$ and 
$(\mu ',g')$ be 
integral diagrams. 
Then $(\mu ,g)$ and $(\mu ',g')$ are
${\mathcal R}^+$-{\it equivalent} if there exist a diffeomorphism germ
$
\Psi :(\R\times (\R\times\R),0)\rightarrow
(\R\times (\R\times\R),0)
$
of the form
$
\Psi (t,x,y)=(t+\alpha (x,y),\psi (x,y))
$
and 
a diffeomorphism germ
$\Phi :(\R^2,0)\rightarrow (\R^2,0)
$
such that $\Psi\circ (\mu ,g)=(\mu ',g')\circ\Phi .$  We remark that if $(\mu ,g)$ and
$(\mu ',g')$ are ${\mathcal R}^+$-equivalent by the above diffeomorphisms, then we have
$
\mu (u)+\alpha\circ g(u)=\mu '\circ\Phi (u)
$
and
$
\psi\circ g(u)=g'\circ\Phi (u)
$
for any $u\in (\R^2,0).$  Thus the diagram $(\mu +\alpha\circ g,g)$
is strictly equivalent to $(\mu ',g').$
\par
We now define the corresponding equivalence relation among Legendrian unfoldings.
Let $\ell _{(\mu ,f)},\ \ell _{(\mu ',f')}:(\R^2,0)\rightarrow 
(J^1(\R\times\R,\R),z_0)$
be Legendrian unfoldings.  We say that $\ell _{(\mu ,f)}$ and $\ell _{(\mu ',f')}$ 
are $S.P^+$-{\it Legendrian
equivalent 
(respectively, $S.P$-Legendrian equivalent)} if there exist a 
contact diffeomorphism germ
$
K:(J^1(\R\times\R,\R),z_0)\rightarrow
(J^1(\R\times\R^n,\R),z'_0),
$
a diffeomorphism germ
$
\Phi :(\R^2,0)\rightarrow (\R^2,0)
$
and a diffeomorphism germ
$
\Psi :(\R\times (\R\times\R),\Pi (z_0))\rightarrow
(\R\times(\R\times\R),\Pi (z'_0))
$
of the form $\Psi (t,x,y)=(t+\alpha (x,y),\psi (x,y))$ (respectively,
$\Psi (t,x,y)=(t,\psi (x,y))$) 
such that
$
\Pi\circ K=\Psi\circ\Pi
$
and
$
 K\circ{\mathcal L}={\mathcal L}'\circ\Phi.
$
It is clear that if $\ell _{(\mu ,f)}$ and $\ell _{(\mu ',f')}$ are $S.P^+$-Legendrian equivalent
(respectively, $S.P$-Legendrian equivalent), then $(\mu ,\pi\circ f)$ and 
$(\mu ',\pi\circ f')$ are ${\mathcal R}^+$-equivalent (respectively, strictly equivalent).
By Theorem 1.1 in \cite{IK}, the converse is also true for generic $(\mu ,f)$ and
$(\mu ',f').$
The notion of the stability of Legendrian unfoldings with respect to
$S.P^+$-Legendrian equivalence (respectively,
$S.P$-Legendrian equivalence) is analogous to the usual notion of the stability
of Legendrian immersion germs with respect to Legendrian equivalence
(cf. Part $I\! I\! I$ in \cite{Ar1}).
\par

On the other hand, 
we can interpret the above equivalence relation in terms of generating 
families.
For the purpose, we use some notations and results in 
\cite{Ar1}. 
Let $\tilde F,\tilde G:(\R\times (\R\times\R),0)\rightarrow (\R,0)$
be generating families of Legendrian unfoldings of Clairaut type.  We say that 
$\tilde F$ and $\tilde G$ are P-${\mathcal C}^+$-equivalent (respectively, P-${\mathcal C}$-equivalent)
if there exists  a diffeomorphism germ
$\Phi :(\R\times (\R\times\R),0)\rightarrow (\R\times (\R\times\R),0)$
of the form $\Phi (t,x,y)=(t+\alpha (x,y),\phi _1(x,y),\phi _2(x,y))$ (respectively, $\Phi (t,x,y)=
(t,\phi _1(x,y),\phi _2(x,y))$) such that
$\langle F\circ\Phi\rangle _{{\mathcal E}_{(t,x,y)}}=\langle G\rangle _{{\mathcal E}_{(t,x,y)}}$
where $\langle G\rangle _{{\mathcal E}_{(t,x,y)}}$ is the ideal generated by $G$ in 
the local ring of function germs ${\mathcal E}_{(t,x,y)}$ of $(t,x,y)$-variables.
We also say that $\tilde F(t,x,y)$ is ${\mathcal C}^+$ (respectively, ${\mathcal C}$)-versal deformation
of $f=F|\R\times 0$ if 
$$
{\mathcal E}_t=\langle \frac{df}{dt}\rangle _{\R}+\langle f\rangle _{{\mathcal E}_t}+
\langle \frac{\partial F}{\partial x}|\R\times \{(0, 0)\}, 
\frac{\partial F}{\partial y}|\R\times \{(0, 0)\}
,1\rangle _{\R}
$$
(respectively,
$$
{\mathcal E}_t=\langle f\rangle _{{\mathcal E}_t}+
\langle \frac{\partial F}{\partial x_1}|\R\times \{0\},\dots ,\frac{\partial F}{\partial x_n}|\R\times \{0\}
,1\rangle _{\R}).
$$
\par
By the similar arguments like as those of Theorems 20.8 and 21.4 in \cite{Ar1}, we can show the following :

\begin{theorem} 
Let 
$\tilde F,\tilde G :(\R\times (\R\times\R),0)\rightarrow (\R,0)$
be generating families of Legendrian unfoldings of Clairaut type $\Phi _F, \Phi _G$
respectively.  Then 
\par\noindent
{\rm (1)} $\Phi _F$ and $\Phi _G$ are S.P$^+$ (respectively, S.P)-Legendrian equivalent
if and only if $\tilde F$ and $\tilde G$ are P-${\mathcal C}^+$ (respectively, ${\mathcal C}$)-equivalent.
\par\noindent
{\rm (2)} $\Phi _F$ is S.P$^+$ (respectively, S.P)-Legendrian stable if and only if
$\tilde F$ is a P-${\mathcal C}^+$ (respectively, ${\mathcal C}$)-versal deformation of 
$f=F|\R\times \{0\}.$
\end{theorem}

The following theorem is a corollary of Damon's general versality theorem in 
\cite{Damon}.

\begin{theorem} 
Let $\tilde F,\tilde G :(\R\times (\R\times\R),0)\rightarrow (\R,0)$
be generating families of Legendrian unfoldings of Clairaut type such that
$\Phi _F, \Phi _G$ are S.P$^+$ (respectively, S.P)-Legendrian stable.
Then $\Phi _F, \Phi _G$ are S.P$^+$ (respectively, S.P)-Legendrian
equivalent if and only if 
$f=F|\R\times \{0\},\ g=G|\R\times \{0\}$ are ${\mathcal C}$-equivalent (i.e. 
$\langle f\rangle _{{\mathcal E}_t} =\langle g\rangle _{{\mathcal E}_t}$).
\end{theorem}

Then the classification theory of function germs by the ${\mathcal C}$-equivalence is
quite useful for our purpose.  For each function germ $f:(\R,0)\rightarrow (\R,0),$
we set
\begin{eqnarray*}
{\mathcal C}{\rm -cod}\, (f)&=&{\rm dim}_{\R}{\mathcal E}_t/\langle f\rangle _{{\mathcal E}_t}, \\
{\mathcal C}^+{\rm -cod}\, (f)&=&{\rm dim}_{\R}{\mathcal E}_t/\langle f\rangle _{{\mathcal E}_t}+
\langle \frac{df}{dt}\rangle _{\R}, \\
{\mathcal K}{\rm -cod}\, (f)&=&{\rm dim}_{\R}{\mathcal E}_t/\langle f\rangle _{{\mathcal E}_t}+
\langle \frac{df}{dt}\rangle _{{\mathcal E}_t}.
\end{eqnarray*}
Then we have the following well-known classification (cf. \cite{martinet}).

\begin{lemma} 
Let $f:(\R,0)\rightarrow (\R,0)$ be a function germ with
${\mathcal K}{\rm -cod}\, (f)<\infty .$  Then $f$ is ${\mathcal C}$-equivalent to the map germ
$t^{\ell +1}$ for some $\ell\in {\mathbb N}.$
\end{lemma}

By the direct calculation, we have
\begin{eqnarray*}
{\mathcal C}{\rm -cod}\, (t^{\ell +1})&=&\ell +1, \\
{\mathcal C}^+{\rm -cod}\, (t^{\ell +1})&=&\ell .
\end{eqnarray*}
Thus we can easily determine ${\mathcal C}$ (respectively, ${\mathcal C}^+$)-versal
deformations of the above germs by the usual method as follows :
\par\noindent
The ${\mathcal C}$-versal deformation :
$$
t^{\ell +1}+\sum_{i=0}^\ell u_{i+1}t^i.
$$
The ${\mathcal C}^+$-versal deformation :
$$
t^{\ell +1}+\sum_{i=0}^{\ell -1}u_{i+1}t^i.
$$
\par
We now ready to give a proof of Theorem \ref{Clairaut}.
\par\noindent
{\it Proof of Theorem \ref{Clairaut} \/}\
 Let $(\mu ,f)$ be a  
Clairaut type equation and the corresponding Legendrian unfolding
$\ell _{(\mu ,f)}$.
By the previous arguments, we have a function germ
$F:(\R\times\R ,0)\lon (\R,0)$ such that  ${\rm Image}\, j^1_1F={\rm Image }\,\ell _{(\mu ,f)}.$
Therefore we consider generic property of $F(t,x).$
By definition $j^1_1F$ is an immersion germ if and only if 
$$
\left(\frac{\partial F}{\partial t}(0),
\frac{\partial ^2F}{\partial t\partial x}(0)\right)\not= 0.
$$
Under this condition, we have the characterization of the fold point and the cusp point of
$\pi\circ j^1_1F$ as follows (cf., \cite{Gibson}\cite{GG}):
\par\noindent
(A) $\pi\circ j^1_1F$ is the fold germ if and only if 
$\dfrac{\partial F}{\partial t}(0) = 0$ 
and 
$\dfrac{\partial ^2F}{\partial t^2}(0) \not= 0. $
\par\noindent
(B) $\pi\circ j^1_1F$ is the cusp germ if and only if 
$$
\dfrac{\partial F}{\partial t}(0) = \dfrac{\partial ^2F}{\partial t^2}(0) = 0  
{\mbox{\rm \ and \ }} 
\dfrac{\partial ^2F}{\partial t\partial x}(0)
\dfrac{\partial ^3F}{\partial t^3}(0)\not=0.
$$
\par
When $j^1_1F$ is not a immersion germ, we have the following 
characterization of the cross cap:
\par\noindent
(C) $j^1_1F$ is a cross cap germ if and only if
$$
\dfrac{\partial F}{\partial t}(0) = 
\dfrac{\partial ^2F}{\partial t\partial x}(0) = 0 
{\mbox{\rm \ and \ }}  
\dfrac{\partial ^3F}{\partial t\partial x^2}(0)
\dfrac{\partial ^2F}{\partial t^2}(0)\not=0.
$$
\par
In the first place, we give normal forms under the assumption that the conditions (a),(b),(c).
Let assume that the condition (C) holds.
In this case the function germ has the following form:
$$
F(t,x)=at^2+bx^2+ctx^2+h(t,x),
$$
where $a\not=0,c\not=0$ and $h(0,0)=0.$
Since $F(t,0)=at^2+h(t,0)$ is ${\mathcal C}$-equivalent to $t^2,$
$F(t,x)$ is $P$-${\mathcal C}$-equivalent to a deformation of $t^2.$
By the previous arguments, the ${\mathcal C}$-versal deformation of $t^2$
is $t^2+v_1t+v_2.$
Therefore, $F(t,x)$ is $P$-${\mathcal C}$-equivalent to
the function germ of the form:
$$
G(t,x)=t^2+t\phi _1(x)+\phi _2(x).
$$
Since $j^1_1G$ is also a cross cap germ, we have
$$
\phi _1(x)=\alpha x^2+{\rm higher\  order\  term},
$$
with $\alpha \not=0.$
By a local diffeomorphim of the variable $x,$
we have
$\phi _1(x)=x^2.$
This means that $F(t,x)$ is $P$-${\mathcal C}$-equivalent to
the germ of the form $t^2+tx^2+\phi (x).$
Hence, we might put that
$F(t,x)=t^2+tx^2+\phi (x).$
In this case.
$$
j^1F(t,x)=(t,x,t^2+tx^2+\phi (x),2t+x^2,2tx+\phi '(x)).
$$
The corresponding integral diagram is
$$
\mu (u_1,u_2)=u_1,\quad g(u_1,u_2)=(u_2,u_1^2+u_1u_2^2+\phi (u_2)).
$$
On the $(x,y)$-plane, we have a diffeomorphism germ
$\Psi :(\R^2,0)\lon (\R^2,0)$ defined by
$\Psi(x, y) = (x, \dfrac{1}{4}(y + \dfrac{1}{4}x^4 - \phi(x)))$. 
Then 
we have $\Psi\circ g(u_1,u_2)=(u_2, \dfrac{1}{4}(u_1 + \dfrac{1}{2}u_2^2)^2).$
This is the normal form (4) in Theorem \ref{Clairaut}, 
after setting $(u, v) = (u_2, u_1 + \dfrac{1}{2}u_2^2)$. 
\par
For the case (A), we can apply  almost the same arguments as the above
and get the normal form of (2) in Theorem \ref{Clairaut}.
\par
For the case (B), the situation is a rather different.
In this case the function $F(t,0)$ is ${\mathcal C}$-equivalent to 
$t^3.$
The ${\mathcal C}$-versal deformation of $t^3$
is $t^3+v_1t^2+v_2t+v_3,$
then the above arguments cannot work in this case.
However, the ${\mathcal C}^+$-versal deformation of $t^3$ is
$t^3+v_1t+v_2.$
Thus we can apply almost the same arguments as the above and
the corrersponding integral diagram is ${\mathcal R}^+$-equvialent to
$$
\mu (u_1,u_2) = u_2,\quad
g(u_1,u_2) = (u_1, u^3_2+u_1u_2).
$$
This means that the diagram is strictly equivalent to the normal form (3) in 
Theorem \ref{Clairaut}.
\par
We now show that the set of function $F(t,x)$ satisfying the conditions (A), (B), (C) or (R) at any point are 
generic in the space of all functions (equiped with the Whitney $C^\infty$-topology).
Here the condition (R) is
that
$\displaystyle{\frac{\partial F}{\partial t}(0)\not= 0.}$
Let $J^3(2,1)$ be the set of $3$-jets of function germs $h:(\R^2 ,0)\lon (\R,0).$
We consider the following two algebraic subset of $J^3(2,1):$
$$
\Sigma _1=\left\{j^3h(0)\ |\ \frac{\partial h}{\partial t}(0)=\frac{\partial^2 h}{\partial t\partial x}(0)=
\frac{\partial ^2h}{\partial t^2}(0)\frac{\partial ^3h}{\partial t\partial x^2}(0)=0\ \right\},
$$ 
$$
\Sigma _2=\left\{j^3h(0)\ |\ \frac{\partial h}{\partial t}(0)=\frac{\partial^2 h}{\partial t^2}(0)=
\frac{\partial ^2h}{\partial t\partial x}(0)\frac{\partial ^3h}{\partial t^3}(0)=0\ \right\}.
$$
We consider the union 
$W =\Sigma _1\cup \Sigma _2,$
then it is also an algebraic subset of $J^3(2,1).$
We can stratify the algebraic set $W$ by submanifolds whose codimensions are at least $3.$
By Thom's jet transversality theorem, $j^3F(\R^2)\cap (\R^2\times\R\times W)=\emptyset $
for a generic function $F(t,x).$
We can easily show that the conditions (A),(B),(C) or (R) are satisfied for such a function $F(t,x).$
This completes the proof of Theorem \ref{Clairaut} and in particular we have 
Therem \ref{gen3}.
\enD

{\small 
\begin{flushleft}
A.A. Davydov; 
Department of Mathematics, Vladimir State University, Gorkii street 87, 
600026 Vladimir, Russia. \\
e-mail : davydov-m2@vpti.vladimir.su \\

\smallskip

G. Ishikawa; 
Department of Mathematics, Hokkaido University, Sapporo 060-0810, Japan. \\
e-mail : ishikawa@math.sci.hokudai.ac.jp \\

\smallskip

S. Izumiya; 
Department of Mathematics, Hokkaido University,  Sapporo 060-0810, Japan. \\
e-mail : izumiya@math.sci.hokudai.ac.jp \\

\smallskip

W.-Z. Sun; 
Department of Mathematics, North East Normal University, 
Chang Chun 130024, P.R. China. \\
e-mail : wzsun@nenu.edu.cn
\end{flushleft}
}

\end{document}